\documentclass[11pt]{amsart}
\numberwithin{equation}{section}
\newtheorem{theorem}{Theorem}
\newtheorem{definition}[theorem]{Definition}
\newtheorem{proposition}[theorem]{Proposition}

\numberwithin{theorem}{section}

\newtheorem{remark}[theorem]{Remark}
\def\R{\bf R}

\def\al{\aligned}
\def\eal{\endaligned}
\def\M{{\bf M}}
\def\be{\begin{equation}}
\def\ee{\end{equation}}
\def\lab{\label}

\def\e{\epsilon}

\def\R{{\bf R}}
\def\M{{\bf M}}

\def\al{\aligned}

\def\p{\partial}
\def\d{\nabla}
\def\lam{\lambda}
\numberwithin{equation}{section}

\begin{document}

\tracingpages 1
\title[minimizer]{\bf minimizers
of the sharp Log entropy on manifolds with non-negative Ricci curvature and flatness
}
\author{Qi S. Zhang}
\address{
Department of
Mathematics, University of California, Riverside, CA 92521, USA}
\date{2017/07/25}

\begin{abstract}

Consider the scaling invariant, sharp log entropy (functional) introduced by Weissler \cite{W:1} on noncompact manifolds with nonnegative Ricci curvature.  It can also be regarded as
a sharpened version of Perelman's W entropy \cite{P:1} in the stationary case. We prove that it has a minimizer if and only if the manifold is isometric to $\R^n$.

Using this result, it is proven that a class of noncompact manifolds with nonnegative Ricci curvature is isometric to $\R^n$. Comparing with the well known flatness results in \cite{An:1}, \cite{Ba:1} and \cite{BKN:1}  on asymptotically flat manifolds and  asymptotically locally Euclidean (ALE) manifolds, their
decay or integral condition on the curvature tensor is replaced by the condition that the metric converges to the Euclidean one in $C^1$ sense at infinity. No second order condition on the metric is needed.

\end{abstract}
\maketitle
\tableofcontents
\section{Statement of result}

Finding extremals of useful functionals, entropies and inequalities is often an useful task in mathematics. Examples include the Sobolev inequality, Perelman's F and W entropies, log Sobolev inequalities, Yamabe functional etc.
In this note we consider the scaling invariant log entropy (functional) introduced by Weissler [We] on noncompact manifolds with nonnegtaive Ricci curvature. It is a
scaling invariant version of the log Sobolev functional originally
introduced by Gross \cite{Gs:1} and Federbush \cite{F:1}. It can also be regarded as
a sharpened version of Perelman's W entropy for the Ricci flow in the stationary case.

\begin{definition}

(a). (After Weissler) Let $\M$ be a Riemannian $n$ manifold. The scaling invariant log functional (entropy) is
\be
\al
L=L(v,
g) &= - \int_\M v^2 \ln v^2 dg +  \frac{n}{2} \ln
\left(\int_\M 4 |\nabla v |^2  dg  \right) + s_n \\
&\equiv - N(v) +   \frac{n}{2} \ln (F(v)) + s_n.
\eal
\ee Here
 $s_n=-\frac{n}{2} \ln (2 \pi n)
- \frac{n}{2}$ and $v \in W^{1, 2} (\M)$ with $\Vert v \Vert_2 =1$;

(b). The infimum of the log Sobolev functional is denoted by
\[
\lam = \lam(g, \M) =
\inf \{ L(v, g) \, | \,  v \in W^{1, 2}_0(\M, g), \quad \Vert v \Vert_{L^2(\M)}=1 \}.
\]

(c). The infimum of the log Sobolev
functional at infinity, in case $\M$ is noncompact, is
\[
\lam_\infty= \lam_\infty(g, \M) = \lim_{r \to \infty} \lam(g, \M-B(x_0, r))
\]where $x_0$ is a reference point in $M$.

\end{definition}

  When $\M= \bf R^n$, then $L(v, g)$ is introduced by Weissler \cite{W:1}.  Some existence results for minimizers of a functional similar to $L(v, g)$ (modified with a scalar curvature term) were proven in (\cite{Z14:1}). Since the functional $L(v, g)$ is scaling invariant, the proof involves an approximation and blow up analysis which is not needed for the usual W functional.   Some applications were given on breathers and irreversibility of world sheets. Existence of minimizers for the W functional on compact manifolds was proven by \cite{R:1}.
  See \cite{Z12:1} for existence and nonexistence results of minimizers of the W functional on noncompact manifolds.

\begin{remark} In the above definition, one can write $u=v^2$
\be
\lab{defNv0}
N=N(v)=\int_\M v^2 \ln v^2 dg = \int_\M u \ln u dg
\ee is just Boltzmann's entropy; and
\be
\lab{defFv0}
F=F(v) = \int_\M  4 |\nabla v |^2   dg = \int_\M \frac{|\nabla u |^2}{u} dg
\ee is just Perelman's F entropy minus a term involving the scalar curvature.

\end{remark}

The first main result of the note is:

\begin{theorem}
\lab{thmjixiao}
 Let $\M$ be a noncompact, complete manifold of dimension $n \ge 3$  such that $Ricci \ge 0$.
The functional $L$ has a minimizer if and only if $\M$ is isometric to $\R^n$.
\end{theorem}

Normally one would believe that a minimizer exists for many "nice" manifolds with nonnegative Ricci curvature. The theorem unexpectedly shows that the only "nice" one is $\R^n$. The proof is given in the next section. An application on flatness of some noncompact manifolds
is given in Section 3.

\section{proof}

\noindent Proof of Theorem \ref{thmjixiao}.

It is well known that the Gaussian functions are minimizers for $L$ in $\R^n$. So one only needs to prove that if $L$ on $\M$ has a minimizer then $\M$ is isometric to $\R^n$.

Let $v$ be a minimizer of $L$ on $\M$ such that $Ric \ge 0$. Then we solve the backward heat equation with final value at time $t=1$.
\be
\lab{bhe}
\al
\begin{cases}
&\Delta u + \p_t u =0 , \quad \text{on} \quad \M \times [0, 1];\\
&u(x, 1) = v^2(x).
\end{cases}
\eal
\ee  The time interval $[0, 1]$ is chosen for convenience. Any finite interval also works.

The minimizer $v$ satisfies the equation (cf. Theorem 1.9 \cite{Z14:1}):
\be
\lab{maineq}
 \frac{n}{2} \frac{4 \Delta v }{ \int ( 4 | \nabla v |^2 ) dg }
+  2 v \ln v +   \left( \lam(g,  \M)+   \frac{n}{2} -  \frac{n}{2}  \ln  \int ( 4 | \nabla v |^2  ) dg
- s_n \right)  v = 0.
\ee Comparing with that paper, the scalar curvature term is dropped here.

According to Lemma 2.3 in \cite{Z12:1}, the function $v$ decays (quadratic ) exponentially near infinity. Hence, by the standard upper bound on heat kernel \cite{LY:1}, we know that $u=u(\cdot, t)$ also decays exponentially for each fixed $t$.  Hence we can use integration by parts to deduce
\[
\p_t N = \p_t \int_\M u \ln u dg = \int_{\M} \p_t u \ln u dg + \int_{\M} \p_t u dg
=- \int_{\M} (\Delta u)  \, \ln u dg = \int_{\M} \frac{|\nabla u|^2}{u} dg = F.
\] Here we should  mention that in the paper \cite{Z12:1} ans \cite{Z14:1}, we assumed that the curvature tensor and/or its gradients are bounded. They are used to deal with the extra scalar curvature term. In this paper, since there is no scalar curvature term, we no longer need these assumptions. All we need are Li-Yau type bound for heat kernel and Hamilton type gradient bound
\cite{Ha:1}. Both hold when $Ric \ge 0$.  Another difference from these two papers is that 
it was assumed that the volume is non-collapsed there.
 But that is not really needed here since
the coefficient in front the exponential term can be made to depend on the volume of the unit ball at one reference point, say $x_0$. i.e.
\[
v(x) \le C_1 e^{- C_2 d^2(x, x_0)}
\]where $C_1$ and $C_2$ are positive constants with $C_1$ depending on $|B(x_0, 1)|$.
Below we give the details of the integration by parts.

 The backward heat kernel $p(x, t; y, 1)$ with $t<1$ is the same as the (usual forward) heat kernel $G=G(y, 1; x, t)$. Since the manifold has $Ricci \ge 0$, Li-Yau's \cite{LY:1}  heat kernel bound implies
\be
\lab{lyjie}
\frac{c^{-1}_3}{|B(y, \sqrt{1-t})|} e^{-c^{-1}_4 d^2(x, y)/(1-t)} \le p(x, t; y, 1) \le \frac{c_3}{|B(y, \sqrt{1-t})|} e^{-c_4 d^2(x, y)/(1-t)}.
\ee Therefore
\[
\al
u(x, t) &= \int p(x, t; y, 1) v^2(y) dy\\ 
&\le  \int \frac{c_3 C^2_1 }{|B(y, \sqrt{1-t})|} e^{-c_4 d^2(x, y)/(1-t)}  e^{- 2 C_2 d^2(y , x_0)} dy\\
 &\le  \int \frac{c_3 C^2_1 }{|B(y, \sqrt{1-t})|} e^{-c_4 d^2(x, y)/[2(1-t)]} 
 e^{-c_4 d^2(x, y)/[2(1-t)]}  e^{- 2 C_2 d^2(y , x_0)} dy\\
 &\le c_5 e^{-c_6 d^2(x , x_0)}, \qquad \forall t \in [0, 1).
\eal
\]Likewise
\be
\lab{dujie}
|\d u(x, t) | \le \frac{c_5}{\sqrt{1-t}} e^{-c_6 d^2(x , x_0)}, \qquad \forall t \in [0, 1).
\ee Next we can use the local gradient bound in \cite{SZ:1}  backwardly for the backward heat equation, which is the same as for the heat equation, to deduce
\[
\frac{|\d u|^2}{u}(x, t) \le u(x, t)  \frac{C}{1-t} \ln^2 \frac{A}{u(x, t)}
\]where $A=\sup_{B(x, 2 \sqrt{1-t}) \times [t, (1+t)/2]} u$. Using (\ref{lyjie}) and volume doubling condition, direct computation shows that
\be
\lab{du2/ujie}
\frac{|\d u|^2}{u}(x, t) \le \frac{c_7}{(1-t)^3}  e^{-c_6 d^2(x , x_0)}, \quad t \in [0, 1).
\ee Here $c_7$ may depend on $|B(x_0, 1)|$. Also the dependence on $(1-t)^{-3}$ can be improved but there is no need to do it here.
Now for any large $r>0$, for each fixed $t \in [0, 1)$, integration by parts shows 
\[
\int_{B(x_0, r)} (\Delta \ln u) u dg = -\int_{B(x_0, r)} \frac{|\d u|^2}{u} dg + 
\int_{\partial B(x_0, r)} \frac{\partial u}{\partial n} dS.
\]Since $|\partial B(x_0, r)| \le C_n r^{n-1}$ by Bishop-Gromov volume comparison, letting $r \to \infty$, we can use (\ref{dujie}) and (\ref{du2/ujie}) to conclude that
\[
\int_\M (\Delta \ln u) u dg = -\int_\M \frac{|\d u|^2}{u} dg,
\]justifying the integration by parts.

Similar to Hamilton's calculation \cite{Ha:1} for the forward heat equation case, we also have
\[
(\Delta + \p_t) (\frac{|\d u|^2}{u}) = 2 | Hess \, \ln u|^2 u + 2 Ric(\d u , \d u)/u.
\]Hence
\[
\p_t F = 2 \int_\M  |Hess \, \ln u|^2 u dg + 2 \int_{\M} Ric(\d u , \d u)/u dg.
\]Again the integration is justified due the exponential decay of $u$.
Therefore
\be
\lab{dtL}
\al
\p_t L &= \left[-F^2 + n  \int_\M  |Hess \, \ln u|^2 u dg + n \int_{\M} Ric(\d u , \d u)/u dg \right] F^{-1}\\
&=F^{-1}\left[ -F^2 + \int_{\M} (\Delta \ln u)^2 u dg \right] +
 n F^{-1} \int_\M  |Hess \, \ln u - \frac{1}{n} (\Delta \ln u )g |^2 u dg \\
&\qquad + n F^{-1} \int_{\M} Ric(\d u , \d u)/u dg.
\eal
\ee

Using integration by parts and Cauchy-Schwarz inequality,
\be
\lab{f2}
F^2 = \left( \int_{\M} \frac{|\nabla u|^2}{u} dg \right)^2 =
\left( \int_{\M} (\Delta \ln u) u   dg \right)^2 \le \int_{\M} (\Delta \ln u)^2 u dg \, \int u dg=\int_{\M} (\Delta \ln u)^2 u dg.
\ee The equality is reached only if
\[
(\Delta \ln u ) \sqrt{u} = C \sqrt{ u}
\]i.e. $\Delta \ln u = C$, since $u>0$.
 As explained above,  integration by parts can be justified by modifying the proof of  Corollary 4.1 in \cite{Z12:1}.

 Hence all three terms on the right hand side of (\ref{dtL}) are non-negative.
From definition, we know $L(v, g) = L(\sqrt{u(\cdot, 1)}, g)= \lam$ and $L(\sqrt{u(\cdot, 0)}, g) \ge \lam$.
Hence
\[
0 \ge L(\sqrt{u(\cdot, 1)}, g) - L(\sqrt{u(\cdot, 0)}, g) = \int^1_0 \p_t L (\sqrt{u}, g) dt.
\]Substituting (\ref{dtL}) to the right hand side, we find that
\[
\Delta \ln u(x, t) = C(t),
\]where $C(t)$ is a function of $t$ only;
\[
Hess \, \ln u - \frac{1}{n} (\Delta \ln u )g =0;
\]
\be
\lab{ricdu}
 Ric(\d u , \d u)=0.
 \ee

Since $u(\cdot, t)$ decays exponentially, we see that $C(t)$ is a non zero constant.
The reason is that if
\[
\Delta \ln u =0
\] then $\ln u$ is harmonic.  But $u=v^2$ and $v$, as a solution to equation (\ref{maineq}),  is bounded from above and nonnegative (cf \cite{Z12:1}). Hence $\ln u < C$.
Yau's Liouville theorem  implies $\ln u$ is a constant, which is impossible.

From these equalities, it is known that $\M$ is isometric to $\R^n$. See Tashiro \cite{T:1} Theorem 2 (I, B). Also Naber \cite{N:1} shows, by a different method, that a Ricci flat shrinking gradient soliton is $\R^n$.
 However, one can not assume that  $\M$ is  Ricci flat yet, since (\ref{ricdu}) holds only in the $\d u$ direction. So $\M$ is not yet a gradient Ricci soliton to begin with. Here we give a very short proof for completeness.
We work with $t=0$. Then,
by considering $f=\ln u$ or $f=-\ln u$ we can assume
\be
\lab{hessf}
Hess \, f =  \lam g,
\ee where $\lam$ is a positive constant. Fix any point $p \in \M$ and pick a point $x \in \M$.
Let $r=r(s)$ be a minimal geodesic connecting $p$ and $x$, parameterized by arc length. Then the definition of Hessian and (\ref{hessf}) tell us
\[
\lam = (Hess \, f) (\p_r, \p_r) = \d_{\p_r} (\d_{\p_r} f) - (\d_{\p_r} \p_r) f= \frac{d^2}{ds^2}
f(r(s)).
\]Therefore, for $r=d(x, p)$,
\[
f(x) = \frac{\lam}{2} r^2 + r \, \frac{d}{ds}
f(r(s)) \, |_{s=0} \ + f(p).
\]This shows that $f$ must have a global minimum.  Choose $p$ to be a minimal point, then
\[
f(x) = \frac{\lam}{2} r^2  + f(p).
\]Note that the smoothness of $f$ implies that $r^2$ is smooth.
Substituting this to (\ref{hessf}) and taking trace, we see that
\[
\frac{\lam}{2} \Delta r^2= \lam n.
\]Hence, the following holds in the classical sense for all $r>0$:
\[
\Delta r = \frac{n-1}{r}.
\]Let $w$ be the volume element in a spherical coordinate centered at $p$. Then
\[
\p_r \ln w = \Delta r = \frac{n-1}{r}.
\]Therefore
\[
\p_r (w/r^{n-1}) = 0.
\]This shows, since $\M$ is a smooth manifold,  $w= w_n r^{n-1}$ where $w_n$ is the volume of standard unit sphere in $\R^n$. Hence the metric is Euclidean by the equality case of the Bishop-Gromov volume comparison theorem.
\qed

\section{flatness of some manifolds with $Ric \ge 0$.}

Next we apply the theorem to the study of flatness of manifolds with $Ric \ge 0$.
Let us recall the definition of asymptotically flat manifolds (cf p64 \cite{LP:1}).

\begin{definition}
\lab{defaf}

A complete, noncompact Riemannian manifold $M$ is called Asymptotically Flat of order $\tau$ if
there is a partition $M=M_0 \cup M_\infty$, which satisfies the following properties.

(i). $M_0$ is compact.

(ii). $M_\infty$ is the disjoint union of finitely many components each of which is diffeomorphic to $({\R^n} - B(0, r_0))$ for some $r_0>0$.

(iii). Under the coordinates induced by the diffeomorphism,  the metric $g_{ij}$ satisfies, for $x \in
M_\infty$,
\[
g_{ij}(x) = \delta_{ij}(x) + O(|x|^{-\tau}), \quad \partial_k
g_{ij}(x) =  O(|x|^{-\tau-1}), \quad \partial_k  \partial_l
g_{ij}(x) =  O(|x|^{-\tau-2}).
\]
\end{definition}

\begin{remark}
\lab{re2.1}
For convenience we will equip the compact component $M_0$ with a reference point $0$.
We will also assume  that $M_\infty$ has only one connected
component. This assumption does not reduce any generality
\end{remark}

 These class of manifolds are quite useful in general relativity and differential geometry. Ricci flat AF manifolds are often the blow up limits in many situations. If one can show these manifolds are $\R^n$, then one usually can prove some useful results by the method of contradiction.
 In \cite{An:1}, \cite{Ba:1} and \cite{BKN:1},   these authors showed that Ricci flat (or $Ricc \ge 0$) AF manifolds  are isometric to
 $\R^n$ if the curvature tensor is in $L^{n/2}$ or it decays faster than inverse square of the distance function. Note there are various definitions of AF manifolds. In the definition used in this paper, AF manifolds are special cases of asymptotically locally Euclidean (ALE) manifolds (cf. \cite{BKN:1}). AF manifolds are ALE manifolds which are simply connected at infinity.
 Related questions and results on flatness of manifolds with nonnegativity of certain curvatures and faster than quadratic curvature decay can be found in \cite{BGS:1} (p58-59), \cite{GW:1}, \cite{K:1} and
 \cite{KS:1}. See also \cite{MSY:1}, \cite{Ni:1} and \cite{NT:1} for results on K\"ahler manifolds.
 
 Here we replace the AF condition by a much weaker asymptotic condition and remove the curvature condition. 

\begin{definition}
\lab{defaf}

A complete, noncompact Riemannian manifold $\M$ is called $C^1$ Asymptotically Euclidean (C1AE) if
there is a partition $\M=M_0 \cup M_\infty$, which satisfies the following properties.

(i). $M_0$ is compact.

(ii). $\M_\infty$ is the disjoint union of finitely many components each of which is diffeomorphic to $({\R^n} - B(0, r_0))$ for some $r_0>0$.

(iii). Under the coordinates induced by the diffeomorphism,  the metric $g_{ij}$ satisfies, for $x \in
M_\infty$,
\[
g_{ij}(x) = \delta_{ij}(x) + o(1), \quad \partial_k
g_{ij}(x) =  o(1).
\]
\end{definition}

Here $o(1)$ means a quantity that goes to $0$ as $|x| \to \infty$.  Observe that there is no decay assumption on the second order derivatives of the metric and hence no decay assumption on curvature. The result of this section is:

\begin{theorem}
A $C^1$ asymptotically Euclidean manifold with $Ric \ge 0$ is isometric to $\R^n$.
\proof
\end{theorem}

  If $\lam=\lam(g, \M)=0$, then it is known from Corollary 1.6 \cite{BCL:1}, that
$\M$ is isometric to $\R^n$.

 So we can assume $\lam<0$. 

We will show that the functional $L$ has a minimizer. Then Theorem \ref{thmjixiao} will imply that $\M=\R^n$ and hence $\lam=0$, reaching a contradiction. So $\M=\R^n$ to begin with.

According to Theorem 1.9 in \cite{Z14:1}, if we can show that
\be
\lab{lam<lam8}
-\infty < \lam < \lam_\infty=0
\ee then a minimizer exists. So we are left to prove (\ref{lam<lam8}).
We mention that the log functional in \cite{Z14:1} has an extra scalar curvature term comparing with the current one. However, since the scalar curvature is non-negative, the same conclusion holds and the proof is the same.

First we prove the following

Claim. {\it  Let $(\M, g)$ be an C1AE manifold of dimension $n \ge
3$.

(a). Then there exists a constant $A>0$,  such that
\be
\lab{alesob} \left( \int_{\M} v^{2n/(n-2)} dg \right)^{(n-2)/n} \le
A \int_{\M} ( 4 |\nabla v |^2 + R v^2 ) dg, \quad \forall  v \in
W^{1, 2}(\M, g);
\ee
moreover $\lam(g)$ is bounded from below i.e.
\be \lab{alelogsob}
 \int_{\M} v^2 \ln v^2 dg \le \frac{n}{2} \ln \left( A \int_{\M} ( 4 |\nabla v |^2 + R v^2 ) dg
 \right),
 \ee $\forall v \in W^{1, 2}(\M, g), \Vert v \Vert_{L^2(\M, g)}=1.$

(b). $\lam_\infty(g) \ge 0$.
}

(a).  We just need to prove (\ref{alesob}) since (\ref{alelogsob})
follows from Jensen inequality.  A C1AE manifold has maximum volume growth, namely,
\[
|B(x, r)| \ge C r^n
\]for some positive constant $C$ and all $r>0$. Then it is well known that (\ref{alesob}) holds.

\medskip

Now we prove part (b).

\noindent First we prove the following assertion.

{\it When the radius $r$ is sufficiently large, we have
\be
\lab{assert2.1}
\lam(g,  \M-B(0, r)) \ge \lam(g_E, {\R}^n-J(B(0, r)) + o(1).
\ee Here $J$ is the coordinate map near infinity in the definition
of C1AE manifold; $o(1)$ is a quantity whose absolute value goes to
$0$ when $ r \to \infty$; $g_E$ is the Euclidean metric.}

Pick a function $v \in C^\infty_0(\M-B(0, r))$ with $\Vert v
\Vert_{L^2}=1$. Given any $\e>0$, by definition of C1AE manifolds,
for $x \in \M-B(0, r)$ with $r$ sufficiently large, there are the following relations
\be
\lab{2.100}
(1-\e) dx \le dg(x) =\sqrt{ det g(x)} dx \le (1+\e) dx,
\ee
\be
\lab{2.101}
(1-\e) |\nabla_{\R^n} f | \le |\nabla v | \le (1+\e)
|\nabla_{\R^n} f |
\ee where $f = v \circ J^{-1}$ and $J$ is the coordinate map. Also
$\nabla_{\R^n}$ is the Euclidean gradient.  Hence
\be
\lab{2.102}
\int_{\M} ( 4 |\nabla v |^2 + R v^2 ) dg
\ge (1-\e )^2
 \int_{\R^n} 4 |\nabla_{\R^n} f|^2 \sqrt{ det g(x)} dx
\ee  Write $\sqrt{ det g(x)} = w^2$,  then
\be
\lab{2.103}
\al
&\int_{\R^n} 4 |\nabla_{\R^n} f|^2 \sqrt{ det g(x)} dx = \int_{\R^n} 4 |w \nabla_{\R^n} f|^2  dx\\
&=\int_{\R^n} 4 |\nabla_{\R^n}(w f) |^2 dx - 8 \int_{\R^n} f \d w \d ( w f) dx +
 4 \int f^2 |\d w |^2 dx.
\eal
\ee By definition of $C1AE$ manifolds, we know that $|\d w | \le \eta(r)$ where $\eta =\eta(r)$ is a function going to $0$ as $ r \to \infty$. Hence,  we have
\be
\al
\lab{2.104}
\int_{\R^n} 4 |\nabla_{\R^n} f|^2 \sqrt{ det g(x)} dx &\ge (1-\eta(r)) \int_{\R^n} 4 |\nabla_{\R^n}( f w) |^2 dx - 16 \eta^{-1}(r) \int f^2 |\d w|^2 dx,\\
&\ge (1-\eta(r)) \int_{\R^n} 4 |\nabla_{\R^n}( f w) |^2 dx - 16 \eta(r) \int f^2  dx
\eal
\ee which implies
\be
\lab{fv>1-e}
\int_{M} ( 4 |\nabla v |^2 + R v^2 ) dg  \ge (1-\e)^2 (1-\eta(r)) \int_{\R^n} 4 |\nabla_{\R^n}( f w) |^2 dx - C \eta(r).
\ee
Also
\be
\lab{2.105}
\al
\int_M v^2 \ln v^2 dg &=
 \int_{\R^n} (f w)^2 \ln f^2 dx =\int_{\R^n} (f w)^2 \ln (f w)^2 dx - \int_{\R^n} (f w)^2 \ln w^2 dx \\
 &=\int_{\R^n} (f w)^2 \ln (f w)^2 dx + o(1).
 \eal
\ee This and (\ref{fv>1-e}) imply
 that
\be
\lab{2.106}
L(v, g,  \M-B(0, r)) \ge L(f w, g_E,  {\R}^n-J(B(0, r))) + o(1)- C \eta(r) - n \e.
\ee Since $\Vert f w \Vert_{L^2(\R^n)} = 1$, by taking the infimum of this inequality,
it is easy to see that
\be
\lab{2.107}
\lam(g,  \M-B(0, r)) \ge \lam(g_E, 
{\R}^n-J(B(0, r)) + o(1) - C \eta(r) -n \e.
\ee Since $\e$ is arbitrary,  the assertion is proven.

Using $\lam(g_E,  {\R}^n-J(B(0, r)) \ge \lam(g_E,  {\R}^n) =
0$, we see that
\be
\lab{2.108}
\lam_\infty(g) = \lim_{r \to \infty} \lam(g,  \M-B(0, r)) \ge 0.
\ee This proves part (b) of the claim.

 But then
\[
\lam < 0 = \lam_\infty.
\]This is (\ref{lam<lam8}) and hence the Theorem follows.
 \qed

\begin{remark}  From the proof, it is clear that one only needs the condition $\lam_\infty=0$ to get the result. In the Ricci flat case, this condition holds if there is a compact set $K$ such that $\M-K$ is conformal to a domain in $\R^n$. In this situation the $W^{1, 2}$ Sobolev constant is the same as the Yamabe constant which stays the same under conformal change. So the best constant in the log Sobolev constant is also the same. This shows $\lam_\infty=\lam_{\R^n}=0$.
Thus we obtain a coordinate free result:
\end{remark}

\begin{proposition}
A complete, noncompact Ricci flat $n (\ge 3)$ dimensional manifold which is conformal at infinity to a domain in $\R^n$ is isometric to $\R^n$.
\end{proposition}

{\bf Acknowledgment.}
Part of the paper was written when the author was a visiting professor at Fudan University. He is grateful to Fudan University and the Simons Foundation
for their support. He  also wishes to thank  Prof. Zhu Meng for helpful conversations.

\bigskip

\noindent e-mail:  qizhang@math.ucr.edu


\begin{thebibliography}{00}

\bibitem [An]{An:1}  Anderson, Michael T. {\it Ricci curvature bounds and Einstein metrics on
compact manifolds}. J. AMS 2 (1989), no. 3,
455-490.

\bibitem [BCL]{BCL:1}  Bakry, D.; Concordet, D.; Ledoux, M. {\it Optimal heat kernel bounds under logarithmic Sobolev inequalities.} ESAIM Probab. Statist. 1 (1995/97), 391-407
    
\bibitem [BGS]{BGS:1} Ballmann, Werner; Gromov, Mikhael; Schroeder, Viktor,
{\it Manifolds of nonpositive curvature.}
Progress in Mathematics, 61. Birkhäuser Boston, Inc., Boston, MA, 1985. vi+263 pp.

\bibitem [Ba]{Ba:1}  Bartnik, R. {\it The mass of asymptotically flat manifolds}, CPAM, 24 (1986), 661-693.

\bibitem [BKN]{BKN:1}  Bando, S., Kasue, A., Nakajima, H. {\it On construction of coordinates at
infinity on manifolds with fast curvature decay and maximum volume growth}, Invent. Math., 97 (1989) 313-349.






\bibitem [D]{D:1}Drees, G. {\it 
Asymptotically flat manifolds of nonnegative curvature.}
Differential Geom. Appl. 4 (1994), no. 1, 77-90. 



\bibitem [F]{F:1} Federbush, P. {\it Partially Alternate Derivation of a Result of Nelson }
J. Math.  Physics, Vol. 10, no 1 Jan. 1969, 50-53



\bibitem [Gs]{Gs:1} Gross, Leonard , {\it Logarithmic Sobolev inequalities.}  Amer. J. Math. 97
(1975), no. 4, 1061-1083.

\bibitem [GW]{GW:1} Greene, R. E.; Wu, H. {\it Gap theorems for noncompact Riemannian manifolds.} Duke Math. J. 49 (1982), no. 3, 731-756.


\bibitem [Ha]{Ha:1} Hamilton, Richard, {\it A matrix Harnack extimate for the heat equation}, Comm. Anal. Geom. 1 (1993), 225-243.
    
\bibitem [K]{K:1} Kasue, A. {\it Gap theorems for minimal submanifolds of Euclidean space.} J. Math. Soc. Japan 38 (1986), no. 3, 473-492.
    
\bibitem [KS]{KS:1} Kasue, Atsushi; Sugahara, Kunio, {\it  Gap theorems for certain submanifolds of Euclidean spaces and hyperbolic space forms.} Osaka J. Math. 24 (1987), no. 4, 679-704.

\bibitem[LP]{LP:1} Lee, J. and Parker, T., {\it The Yamabe problem}, Bull. AMS, 17  no.1 (1987), 37-91.


\bibitem[LY]{LY:1}
  Peter Li and S. T. Yau, {\it On the Parabolic Kernel of the
  Sch\"{o}dinger operator,} Acta  math. \textbf{156} (1986), pp 153--201.
  
\bibitem[MSY]{MSY:1} Mok, Ngaiming; Siu, Yum Tong; Yau, Shing Tung, {\it The Poincar\'e-Lelong equation on complete K\"ahler manifolds}. Compositio Math. 44 (1981), no. 1-3, 183-218.

\bibitem[N]{N:1} Naber, A. {\it Noncompact shrinking 4-solitons with nonnegative curvature.}
J. Reine Angew. Math. 645 (2010), 125-153.


\bibitem[Ni]{Ni:1} Ni, Lei, {\it An optimal gap theorem. Invent. Math.} 189 (2012), no. 3, 737-761.
    
\bibitem[NT]{NT:1} Ni, L. Tam, L-F. {\it Poincaré-Lelong equation via the Hodge-Laplace heat equation.} Compos. Math. 149 (2013), no. 11, 1856-1870.

\bibitem[P]{P:1} Perelman, Grisha, {\it The entropy formula for the Ricci flow and its
geometric applications}, Math. ArXiv, math.DG/0211159.

\bibitem[R]{R:1}Rothaus, O. S. {\it Logarithmic Sobolev inequalities and the spectrum of Schrödinger operators.} J. Funct. Anal. 42 (1981), no. 1, 110-120.
    
\bibitem[SZ]{SZ:1}  Souplet, P. and   Zhang, Qi S. {\it Sharp gradient estimate and Yau's Liouville theorem for the heat equation on noncompact manifolds},
Bulletin of the London Mathematical Society 38 (6), 1045-1053

\bibitem[T]{T:1} Tashiro, Yoshihiro, {\it  Complete Riemannian manifolds and some vector fields.} Trans. Amer. Math. Soc. 117 (1965) 251-275.





\bibitem[W]{W:1} Weissler, Fred B. {\it Logarithmic Sobolev inequalities for the
heat-diffusion semigroup.} Trans. Amer. Math. Soc. 237 (1978),
255-269.





\bibitem[Z12]{Z12:1} Zhang, Qi S. {\it Extremal of Log Sobolev inquality and W entropy on noncompact manifolds},
  J. Func. Analysis 263 (2012),  2051-2101.

\bibitem[Z14]{Z14:1} Zhang, Qi S. {\it A no breathers theorem for some noncompact Ricci flows.} Asian J. Math. 18
(2014), no. 4, 727-755.



\end{thebibliography}
\enddocument